\documentclass[a4,12pt]{article}
\usepackage{booktabs}
\usepackage{amsmath}
\usepackage{graphics}
\usepackage{graphicx}
\usepackage{epstopdf}
\usepackage{latexsym}
\usepackage{amsfonts}
\usepackage{float}
\usepackage{subcaption}
\usepackage{stmaryrd}
\usepackage{color}
\numberwithin{equation}{section} \rightmargin 1.5cm \leftmargin
-1.5cm \textheight 23cm \textwidth 16cm \topmargin -1cm
\title{
 Stochastic differential switching game in infinite horizon}
\author{Brahim EL ASRI \thanks{Universit\'e Ibn Zohr, Equipe. Aide à la decision,
ENSA, B.P.  1136, Agadir, Maroc. e-mail: b.elasri@uiz.ac.ma }\,\,\,
\, and \, Sehail MAZID \thanks{Universit\'e Ibn Zohr, Equipe. Aide à la decision,
ENSA, B.P.  1136, Agadir, Maroc. e-mail: sehail.mazid@edu.uiz.ac.ma.} }
\begin{document}
\date{}
\maketitle
\newtheorem{theo}{Theorem}
\newtheorem{problem}{Problem}
\newtheorem{pro}{Proposition}
\newtheorem{cor}{Corollary}
\newtheorem{axiom}{Definition}
\newtheorem{rem}{Remark}
\newtheorem{lem}{Lemma}
\newcommand{\brm}{\begin{rem}}
\newcommand{\erm}{\end{rem}}
\newcommand{\beth}{\begin{theo}}
\newcommand{\eeth}{\end{theo}}
\newcommand{\bl}{\begin{lem}}
\newcommand{\el}{\end{lem}}
\newcommand{\bp}{\begin{pro}}
\newcommand{\ep}{\end{pro}}
\newcommand{\bcor}{\begin{cor}}
\newcommand{\ecor}{\end{cor}}
\newcommand{\be}{\begin{equation}}
\newcommand{\ee}{\end{equation}}
\newcommand{\beq}{\begin{eqnarray*}}
\newcommand{\eeq}{\end{eqnarray*}}
\newcommand{\beqa}{\begin{eqnarray}}
\newcommand{\eeqa}{\end{eqnarray}}
\newcommand{\dg}{\displaystyle \delta}
\newcommand{\cm}{\cal M}
\newcommand{\cF}{{\cal F}}
\newcommand{\cR}{{\cal R}}
\newcommand{\bF}{{\bf F}}
\newcommand{\tg}{\displaystyle \theta}
\newcommand{\w}{\displaystyle \omega}
\newcommand{\W}{\displaystyle \Omega}
\newcommand{\vp}{\displaystyle \varphi}
\newcommand{\ig}[2]{\displaystyle \int_{#1}^{#2}}
\newcommand{\integ}[2]{\displaystyle \int_{#1}^{#2}}
\newcommand{\produit}[2]{\displaystyle \prod_{#1}^{#2}}
\newcommand{\somme}[2]{\displaystyle \sum_{#1}^{#2}}
\newlength{\inter}
\setlength{\inter}{\baselineskip} \setlength{\baselineskip}{7mm}
\newcommand{\no}{\noindent}
\newcommand{\rw}{\rightarrow}
\def \ind{1\!\!1}
\def \R{I\!\!R}
\def \N{I\!\!N}
\def \cadlag {{c\`adl\`ag}~}
\def \esssup {\mbox{ess sup}}

\begin{abstract}
We study a zero-sum stochastic differential switching game in infinite horizon. We prove the existence of the value of the game and characterize it as
the unique viscosity solution of the associated system of quasi-variational inequalities with bilateral obstacles. We also obtain
a verification theorem which provides an optimal strategy of the game. Finally, some numerical examples with two regimes are given.

\end{abstract}

\no{\bf AMS  subject classifications}:  49N70, 49L25, 60H30, 90C39, 93E20
\medskip

\no {$\bf Keywords$}:  stochastic differential game, switching strategies, quasi-variational inequality, value function, viscosity solution

\section {Introduction}
\no

Differential game theory involves multiple persons (also called players or individuals) decision making in the context of dynamic system. The study of the two-person zero-sum differential games could be traced back to the pioneering work by Isaacs \cite{[IR]}, which inspired further research in this area.

In 1989, Fleming and Souganidis \cite{[FS]}, for the first time, studied two-person zerosum stochastic differential games. They adopted the definitions of upper and lower value functions introduced by Elliott and Kalton \cite{[EKa]} and proved that the two value functions satisfy the dynamic programming principle and they are the unique viscosity solutions to the associated Hamilton–-Jacobi–-Bellman–-Isaacs partial differential equations (HJBI PDEs). Their work generalized that of Evans and Souganidis \cite{[ES]} from the deterministic framework to the stochastic one and is now regarded as one of the outstanding results in the field of stochastic differential games. In 1995, Hamad\`ene and Lepeltier \cite{[HL]} introduced the backward stochastic differential equation (BSDE) theory into the study of stochastic differential games, and there exist also many other works using this BSDEs approach, such as Hamad\`ene, Lepeltier, and Peng \cite{[HLP]}. Subsequently, Buckdahn and Li \cite{[BL]} extended the findings presented in \cite{[HL],[HLP]}, and generalized the framework introduced in \cite{[FS]}.

In this paper  we consider the state process of the stochastic differential game, defined as the solution of the following stochastic equation:
\begin{equation*}
\label{az}
\begin{array}{ll}
 X_{s}=x+\integ{t}{s}b(r,X_{r},a_t,b_t)dr+\integ{t}{s}\sigma(r,X_{r},a_t,b_t)dW_{r}\qquad s\geq t,
\end{array}
\end{equation*}
with $X_{t^-} = x$. Here $W$ is a d-dimensional Wiener process, while
\begin{equation*}
a_t=\sum_{m\geq 1}\xi_{m}\ind_{(\tau_{m}\leq t <\tau_{m+1}]}(t)\quad\quad\text{and}\quad\quad b_t=\sum_{n\geq 1}\eta_{n}\ind_{(\rho_{n}\leq t <\rho_{n+1}]}(t) .
\end{equation*}
with the cost functional
\begin{equation}
\label{eq0}
\begin{array}{ll}
\mathbb{E}\bigg[\integ{0}{\infty}e^{-rs}f(X_{s},a_t,b_t)ds-
\sum\limits_{m\geq 1}e^{-r\tau_m} C(\xi_{m-1},\xi_m) +\sum\limits_{n\geq 1}e^{-r\rho_n}\chi(\eta_{n-1},\eta_n)\bigg].
\end{array}
\end{equation}
The first player chooses the control a from a given finite set $\mathcal{I}$ to maximize the payoff (\ref{eq0}), and each actions is related with one cost C, while the second player chooses the control b from  $\mathcal{I}$ to minimize the payoff (\ref{eq0}), and each of his actions is associated with the other cost $\chi$. The zero-sum stochastic differential games
problems we will investigate is to show the upper and lower value functions coincide and the game admits a value.
The Isaacs' system of
quasi-variational inequalities for this switching game is the following:
for any $i,j\in \mathcal{I}$, and $x\in \mathbb{R}^{m}$,
\begin{equation} \label{eq:HJBI0}
\begin{array}{c}
max\Big\{ min\Big[rV^{ij}(x)-\mathcal{A}V^{ij}(x)-f(x,i,j);\qquad\qquad\qquad\qquad\\ \qquad V^{ij}(x)-M^{ij}[V](x)\Big];V^{ij}(x)-N^{ij}[V](x)\Big\}=0,
\end{array}
\end{equation}
and
\begin{equation} \label{eq:HJBI01}
\begin{array}{c}
min\Big\{ max\Big[rV^{ij}(x)-\mathcal{A}V^{ij}(x)-f(x,i,j);\qquad\qquad\qquad\qquad\\ \qquad V^{ij}(x)-N^{ij}[V](x)\Big];V^{ij}(x)-M^{ij}[V](x)\Big\}=0.
\end{array}
\end{equation}
where,
$$M^{ij}[V](x)= \max_{k\ne i}\{V^{kj}(x)-C(i,k)\}, \qquad N^{ij}[V](x)= \min_{l\ne j}\{V^{il}(x)+\chi(j,l)\}.$$


In the finite horizon framework, the switching game have been studied by several authors. The most recent work
discussing this topic includes the papers by Djehiche et al. (2017) (\cite{[DHMZ]}), Tang and Hou (2007) (\cite{[TH]}).

The objective of this work is to establish existence and uniqueness of a continuous viscosity solution of (\ref{eq:HJBI0}) and (\ref{eq:HJBI01}). The  second contribution of our paper is the Verification Theorem which provides an optimal strategy of the game, To the best of our knowledge, this issue have not been addressed in the literature yet.



This paper is organized as follows: in Section 2 we list all the notations, state the full set of assumptions, and define viscosity sub- and supersolutions along with equivalent characterizations. In Section 3, we shall introduce the stochastic differential game problem and give some preliminary results of the lower and the upper value functions of stochastic differential game. In Section 4, by the
dynamic programming principle we prove that the lower and upper value functions of the game satisfy the Isaacs' system of
quasi-variational inequalities in the viscosity solution sense. In Section 5, we show that the solution of Isaacs' system of
quasi-variational inequalities is unique. Further, the upper and the lower value functions coincide and the game admits a value. In Section 6, we present a verification theorem which gives an optimal strategy of the switching game, while finally Section 7 will present some numerical examples.

\section{Assumptions and problem formulation}
\no

Throughout this paper  $m$ and $d$ are two integers. Let $\mathcal{I}=\{1,...,q\}$ the finite set of regimes, and assume the following assumptions:  \medskip

\no\textbf{[H1]}  $b:\mathbb{R}^{m}\times\mathcal{I}\times\mathcal{I}\rightarrow
\mathbb{R}^{m}$ and $\sigma :\mathbb{R}^{m}\times\mathcal{I}\times\mathcal{I}\rightarrow
\mathbb{R}^{m\times d}$ be two continuous functions for which there
exists a constant $C>0$ such that for any
$x,x^{\prime }\in \mathbb{R}^{m}$ and $i,j\in\mathcal{I}$
\begin{equation}
|\sigma(x,i,j)-\sigma(x^{\prime },i,j)|+|b(x,i,j)-b(x^{\prime},i,j)|\leq
C|x-x^{\prime }|.
 \label{eqs}
\end{equation}%
Thus they are also of linear growth. i.e., there exists a constant $C$ such that for any $x\in\mathbb{R}^m$ and $i,j\in\mathcal{I}$
\begin{equation}
|\sigma(x,i,j)|+|b(x,i,j)|\leq C(1+|x|).
\end{equation}
\textbf{[H2]}\; $f:\mathbb{R}^{m}\times\mathcal{I}\times\mathcal{I} \rightarrow \mathbb{R}$
is a continuous function for which there exists a constant $C$ such that for each $i,j\in\mathcal{I}$, $x,x^{\prime }\in \mathbb{R}^{m}$:
 \begin{equation}
 |f(x,i,j)|\leq C(1+|x|) \qquad\text{and}
 \end{equation}
  \begin{equation}
\qquad\qquad|f(x,i,j)-f(x^{\prime},i,j)|\leq
C|x-x^{\prime }|.\qquad\qquad\quad
\end{equation}
\textbf{[H3]}\; For any $i,j\in\mathcal{I}$ The switching costs $C(i,j)$ and $\chi(i,j)$ are constants, and we assume the triangular condition :
\begin{equation}\label{co1}
 C(i,k)< C(i,j)+C(j,k), \quad j\ne i,k.
\end{equation}
\begin{equation}\label{co2}
\chi(i,k)< \chi(i,j)+\chi(j,k), \quad j\ne i,k,
\end{equation}
which means that it is less expensive to switch directly in one step from regime $i$ to $k$ than
in two steps via an intermediate regime $j$. Notice that a switching costs $C(i,j)$ and $\chi(i,j)$ may be negative,
and conditions (\ref{co1}) and (\ref{co2}) for $i=k$ prevents an arbitrage by simply switching back
and forth, i.e.
\begin{equation}\label{co3}
0< C(i,j)+C(j,i).
\end{equation}
\begin{equation}\label{co4}
0< \chi(i,j)+\chi(j,i).
\end{equation}
We set $b^{ij}(.)=b(.,i,j), \sigma^{ij}(.)=\sigma(.,i,j)$ and $f^{ij}(.)=f(.,i,j)$.

We now consider the following  Isaacs' system of
quasi-variational inequalities: for any $i,j\in \mathcal{I}$, and $x\in \mathbb{R}^{m}$,
\begin{equation} \label{eq:HJBI1}
\begin{array}{c}
max\Big\{ min\Big[rV^{ij}(x)-\mathcal{A}V^{ij}(x)-f^{ij}(x);\qquad\qquad\qquad\qquad\\ \qquad V^{ij}(x)-M^{ij}[V](x)\Big];V^{ij}(x)-N^{ij}[V](x)\Big\}=0,
\end{array}
\end{equation}
and
\begin{equation} \label{eq:HJBI2}
\begin{array}{c}
min\Big\{ max\Big[rV^{ij}(x)-\mathcal{A}V^{ij}(x)-f^{ij}(x);\qquad\qquad\qquad\qquad\\ \qquad V^{ij}(x)-N^{ij}[V](x)\Big];V^{ij}(x)-M^{ij}[V](x)\Big\}=0.
\end{array}
\end{equation}
Where\\
(i) $r$ is a positive discount factor and $\mathcal{A}$ is the following infinitesimal generator:
$$\mathcal{A}V^{ij}= \langle b^{ij},\nabla _x V^{ij}\rangle +\cfrac{1}{2}tr[\sigma^{ij}(\sigma^{ij})^*\nabla^2_xV^{ij}].$$
 (ii) For any $x\in \mathbb{R}^{m}$ and $i,j\in\mathcal{I},$
$$M^{ij}[V](x)= \max_{k\ne i}\{V^{kj}(x)-C(i,k)\}, \qquad N^{ij}[V](x)= \min_{l\ne j}\{V^{il}(x)+\chi(j,l)\}.$$

We now define the notions of viscosity solution of the system (\ref{eq:HJBI1}). We can similarly define the notions for (\ref{eq:HJBI2}).
\begin{axiom}
Let $\vec{V}=(V^{kl}(x))_{(k,l)\in\mathcal{I}\times\mathcal{I}}:\R^m\to\R$ such that for any $(i,j)\in\mathcal{I}\times\mathcal{I}$, $V^{ij} $ is continuous, is called:\\
(i) A viscosity supersolution to (\ref{eq:HJBI1}) if for any $i,j\in\mathcal{I}$, for any $x_0\in\R^m$ and any function $\phi^{ij}\in C^{2}(\R^m)$ such that $\phi^{ij}(x_0)=V^{ij}(x_0)$ and $x_0$ is a local maximum
of $\phi^{ij}-V^{ij}$, we have:
\begin{equation}
\begin{array}{lll}
max\Big\{ min\Big[r\phi^{ij}(x_0)-\mathcal{L}\phi^{ij}(x_0)-f^{ij}(x_0),\\ \qquad V^{ij}(x_0)-M^{ij}[V](x_0)\Big],V^{ij}(x_0)-N^{ij}[V](x_0) \Big\}\geq 0.
\end{array}
\end{equation}\\
(ii) A viscosity subsolution to (\ref{eq:HJBI1}) if for any $i,j\in\mathcal{I}$, for any $x_0\in\R^m$ and any function $\phi^{ij}\in C^{2}(\R^m)$ such that $\phi^{ij}(x_0)=V^{ij}(x_0)$ and $x_0$ is a local minimum
of $\phi^{ij}-V^{ij}$, we have:
\begin{equation}
\begin{array}{lll}
max\Big\{ min\Big[r\phi^{ij}(x_0)-\mathcal{L}\phi^{ij}(x_0)-f^{ij}(x_0),\\ \qquad V^{ij}(x_0)-M^{ij}[V](x_0)\Big],V^{ij}(x_0)-N^{ij}[V](x_0) \Big\}\leq 0.
\end{array}
\end{equation}\\
(iii) A viscosity solution if it is both a viscosity supersolution and subsolution.  \;$\Box$
\end{axiom}

There is an equivalent formulation of this definition (see
e.g.\cite{[CIL]}) which we give since it will be
useful later. So firstly we define the notions of superjet and subjet of a continuous function $V$.
\begin{axiom}
Let $V\in C(\R^m)$, $x$ an element of $\R^m$ and finally $\mathbb{S}_m$ the set of $m\times m$ symmetric matrices. We denote by $J^{2,+}V(x)$ (resp. $J^{2,-}V(x))$, the superjets (resp. the subjets) of\, $V$ at x, the set of pairs $(q,X)\in \R^m\times \mathbb{S}_m$ such that:
\begin{equation*}
\begin{array}{ll}
V(y)\leq V(x)+\langle q,y-x\rangle+\frac{1}{2}\langle X(y-x),y-x\rangle+o(|y-x|^2)\\ \\
(\text{resp.}\; V(y)\geq V(x)+\langle q,y-x\rangle+\frac{1}{2}\langle X(y-x),y-x\rangle+o(|y-x|^2)).\Box
\end{array}
\end{equation*}
Note that if $\phi-V$ has a local maximum (resp. minimum) at $x$, then we obviously have:
$$(D_x\phi(x),D^2_{xx}\phi(x))\in J^{2,-}V(x)\;(\text{resp}.\; J^{2,+}V(x)).\Box$$
\end{axiom}

We now give an equivalent definition of a viscosity solution of of the system (\ref{eq:HJBI1}).
\begin{axiom}
Let $ \vec{V}=(V^{kl}(x))_{(k,l)\in\mathcal{I}\times\mathcal{I}}:\R^m\to\R$ such that for any $(i,j)\in\mathcal{I}\times\mathcal{I}$; $V^{ij} $ continuous , is called a viscosity supersolution (resp. a viscosity subsolution) to (\ref{eq:HJBI1}) if for any $i,j\in\mathcal{I}$, for any $x\in\R^m$ and any $(q,X)\in J^{2,-}V^{ij}(x)$ (resp. $J^{2,+}V^{ij}(x)),$
\begin{equation}
\begin{array}{lll}
max\Big\{ min\Big[rV^{ij}-\frac{1}{2}Tr[(\sigma^{ij})^*X\sigma^{ij}]-\langle b^{ij},q\rangle-f^{ij}(x);\qquad\qquad\qquad\qquad\\ \qquad\qquad V^{ij}(x)-M^{ij}[V](x)\Big];V^{ij}(x)-N^{ij}[V](x) \Big\}\geq 0 \quad(resp. \leq 0).
\end{array}
\end{equation}
It is called a viscosity solution if it is both a viscosity subsolution and supersolution. $\Box$
\end{axiom}

As pointed out previously we will show that system (\ref{eq:HJBI1}) and (\ref{eq:HJBI1}) has a unique solution in viscosity sense.
This systems are the deterministic version of the zero-sum switching game which we will describe briefly in the next section.

\section{The zero-sum switching game}
\subsection{Setting of the problem}
\no

 Consider a complete probability space ($\Omega,\mathcal{F},\mathbb{P}$) and a $d$-dimensional standard Wiener process $W=(W_{t})$ defined on it. Let $\mathbb{F} =(\mathcal{F}_{t})_{t\geq 0}$ be the natural filtration
generated by the Wiener process, completed with the $\mathbb{P}$-null sets of  $\mathcal{F}$. We begin by description of the zero-sum switching game.
\begin{axiom}
Assume we have two players I and II who intervene on a system (e.g. the production of energy from several
sources such as oil, cole, hydro-electric, etc.) with the help of switching strategies. An admissible switching process
for I (resp. II) is a sequence $\delta:=(\tau_m,\xi_m)_{m\geq0}$ (resp. $\nu := (\rho_n,\eta_n)_{n\geq0})$ where:\\
(i) $(\tau_{m})_{m}$ (resp. $(\rho_{n})_{n})$, the action times, is a sequence of $\mathcal{F}$-stopping times such that $\tau_m<\tau_{m+1}$ (resp. $\rho_n<\rho_{n+1})$ and $\tau_n\to\infty$ (resp. $\rho_n\to\infty$).
\\(ii) $(\xi_{m})_{m}$ (resp. $(\eta_{n})_{n})$, the actions, is a sequence of $\mathcal{I}$-valued  random variables, where each $\xi_{m}$ (resp. $\eta_{n}$) is $\mathcal{F}_{\tau_{m}}$-measurable (resp. $\mathcal{F}_{\rho_{n}}$-measurable).
\end{axiom}

Next let $(i,j)\in\mathcal{I}\times\mathcal{I}$ be fixed. We say that the admissible switching strategy $\delta:=(\tau_m,\xi_m)_{m\geq0}$ of I (resp. $\nu :=(\rho_n,\eta_n)_{n\geq0})$ (resp. II) belongs to $\mathcal{A}^i$ (resp. $\mathcal{B}^j$) if $\tau_0=0, \xi_0=i$ (resp. $\rho_0=0, \eta_0=j$).

Given an admissible strategy $\delta$ (resp. $\nu$) of I (resp. II), one associates a stochastic process $(a_t)_{t\geq 0}$ (resp.
$(b_t)_{t\geq0})$ which indicates along with time the current mode of I (resp. II) and which is defined by: $\forall t \geq 0$
 \begin{equation}
a_t=\xi_0\ind_{\{\tau_0\}}(t)+\sum_{m\geq 1}\xi_{m-1}\ind_{(\tau_{m-1},\tau_{m}]}(t)\;(\text{resp.}\; b_t=\eta_0\ind_{\{\rho_0\}}(t)+\sum_{n\geq 1}\eta_{n-1}\ind_{(\rho_{n-1},\rho_{n}]}(t)).
\end{equation}

The evolution of
the state of the game is described by the following stochastic equation:\\
\be \left\{
\begin{array}{lllllll}\label{equa}
dX_t=b^{a_tb_t}(X_t)dt+\sigma^{a_tb_t}(X_t)dW_t, \qquad\; t\geq 0,\\
X_0=x
\end{array}
\right. \ee

We assume that for any $x\in\mathbb{R}^m$ and $i,j\in\mathcal{I}$, there exists a unique strong solution to (\ref{equa}). Let now $\delta:=(\tau_m,\xi_m)_{m\geq0}$ (resp. $\nu := (\rho_n,\eta_n)_{n\geq0})$ be an admissible strategy for I (resp. II) which belongs
to $\mathcal{A}^i$
(resp. $\mathcal{B}^j$). The interventions of the players are not free and generate a payoff which is a reward
(resp. cost) for I (resp. II) and whose expression is given by
\begin{equation}
\label{aar}
\begin{array}{ll}
J(x, \delta, \nu)=\mathbb{E}\bigg[\integ{0}{\infty}e^{-rs}f^{a_sb_s}(X^x_{s})ds-
\sum\limits_{m\geq 1}e^{-r\tau_m} C(\xi_{m-1},\xi_m) +\sum\limits_{n\geq 1}e^{-r\rho_n}\chi(\eta_{n-1},\eta_n)\bigg].
\end{array}
\end{equation}
Following Buckdahn and Li \cite{[BL]}, we define nonanticipative strategy as follows.
\begin{axiom}
A nonanticipative strategy for player I  is a mapping
$$\alpha:\cup_{j\in\mathcal{I}}\mathcal{B}^j\rightarrow\mathcal{A}^i$$
such that for any stopping time $\tau$ and any $b$, $b' \in \mathcal{B}^j$, with $b\equiv b'$ on $\llbracket 0,\tau \rrbracket,$ it holds that $ \alpha(b)\equiv \alpha(b')$ on $\llbracket 0,\tau\rrbracket$.
(with the notation $\llbracket 0,\tau \rrbracket= \{(s,\omega)\in[0,T]\times\Omega,\, 0\leq s \leq \tau(\omega)\}).$ \\ A nonanticipative strategy for player
II
$$ \beta:\cup_{i\in\mathcal{I}}\mathcal{A}^i\rightarrow\mathcal{B}^j$$
is defined similarly. The set of all nonanticipative strategies $\alpha$ (resp. $\beta$) for player I
 (resp. II) is denoted by
$ \Gamma^i$ (resp, $\Delta^j)$.\\
We define the upper (resp. lower) value of the game by
\begin{equation}
\overline{V}^{\;ij}(x)=\inf\limits_{\beta\in {\Delta^j}}\sup\limits_{\delta\in\mathcal{A}^i}J(x,\delta,\beta(\delta))
\end{equation}
\big(resp.
\begin{equation}
\underline{V}^{ij}(x)=\sup\limits_{\alpha\in \Gamma^i}\inf\limits_{\nu\in\mathcal{B}^j}J(x,\alpha(\nu),\nu)\big).
\end{equation}
The game has a value if $\underline{V}^{ij}=\overline{V}^{\;ij}.$
\end{axiom}
\subsection{Preliminary results}
\no

 In this section we present some properties of the lower and upper value functions
of our switching game.
\begin{lem}
(see e.g. \cite{[RY]}) The process $X^{x}$ satisfies the following estimate:
There exists a constant $\rho$
such that
\begin{equation}\label{estimat1}
E[|X^{x}_t|]\leq e^{\rho t}(1+|x|)\qquad \forall t \geq 0.
\end{equation}
\end{lem}
\begin{pro}
Under the standing assumptions (\textbf{H1}), (\textbf{H2}) and (\textbf{H3}), there exists some positive constant $\rho$ such that for $r\geq\rho$,  the lower and upper value function are
satisfy a linear growth condition: for all $i,j\in\mathcal{I}$ and $x\in\R^m$ there exists a constant $C$ such that
\begin{equation}
\label{LG}
  |\overline{V}^{\;ij}(x)|,|\underline{V}^{ij}(x)|\leq C(1+|x|).
\end{equation}
\end{pro}
$Proof$: We make the proof only for the lower value function $\overline{V}^{\;ij}$, the other case being
analogous.
for every $\delta\in\mathcal{A}^i$ we consider the particular strategy:  $\tilde{\beta}(\delta)=(\tilde{\rho}_n,\tilde{\eta}_n)$  given by  $\tilde{\rho}_n=\infty, \tilde{\eta_n}=j$ for all $n\geq 1$  we have
\begin{equation}
\overline{V}^{\;ij}(x) \leq \sup\limits_{\delta\in\mathcal{A}^i}\mathbb{E}\bigg[\integ{0}{\infty}e^{-rs}f^{a_sb_s}(X_{s}^x)ds-
\sum\limits_{m\geq 1}e^{-r\tau_m} C(\xi_{m-1},\xi_m)\bigg].
\end{equation}
Then for every $\epsilon > 0$, there
exists  a strategy $\delta\in \mathcal{A}^i$ such that
\begin{equation}\label{eqboun}
\overline{V}^{\;ij}(x) \leq \mathbb{E}\bigg[\integ{0}{\infty}e^{-rs}f^{a_sb_s}(X_{s}^x)ds-
\sum\limits_{m\geq 1}e^{-r\tau_m} C(\xi_{m-1},\xi_m)\bigg]+\epsilon.
\end{equation}
The switching cost is dominated
by
 \begin{equation}\label{domcost}
-\sum_{m=1}^{N}e^{-r\tau_m}C(\xi_{m-1},\xi_m)\leq \max_{k\in\mathcal{I}}(-C(i,k))
\end{equation}
for any $N \geq 1$, which can be proved by induction as in \cite{[LP]}. Indeed, (\ref{domcost}) obviously holds for
 N= 1. Suppose now that it is also verified for some $N-1$ and let us show that it is valid for $N$. When $C(\xi_{N-1},\xi_N)\geq 0$ holds as $r\geq 0$. When $C(\xi_{N-1},\xi_N)< 0$, we have
\begin{equation*}
\begin{array}{ll}
-\sum_{m=1}^{N}e^{-r\tau_m}C(\xi_{m-1},\xi_m)\\ \\ \qquad\qquad\qquad=-\sum_{m=1}^{N-2}e^{-r\tau_m}C(\xi_{m-1},\xi_m)-(e^{-r\tau_{N-1}}C(\xi_{N-2},\xi_{N-1})+e^{-r\tau_{N}}C(\xi_{N-1},\xi_{N}))\\ \\
\qquad\qquad\qquad\leq -\sum_{m=1}^{N-2}e^{-r\tau_m}C(\xi_{m-1},\xi_m)-e^{-r\tau_{N-1}}(C(\xi_{N-2},\xi_{N-1})+C(\xi_{N-1},\xi_{N}))\\ \\
\qquad\qquad\qquad\leq -\sum_{m=1}^{N-2}e^{-r\tau_m}C(\xi_{m-1},\xi_m)-e^{-r\tau_{N-1}}C(\xi_{N-2},\xi_{N})\leq \max_{k\in\mathcal{I}}(-C(i,k)).
\end{array}
\end{equation*}

Combining (\ref{eqboun}) and (\ref{domcost}) gives us
\begin{equation}
\overline{V}^{\;ij}(x) \leq \mathbb{E}\bigg[\integ{0}{\infty}e^{-rs}f^{a_sb_s}(X_{s}^x)ds+\max_{k\in\mathcal{I}}(-C(i,k))\bigg]+\epsilon.
\end{equation}
Now, from estimate (\ref{estimat1}) and the polynomial growth condition of $f^{ij}$  in (\textbf{H2}), there exists $C$ and $\rho$ such that
\begin{equation}
\overline{V}^{\;ij}(x) \leq \frac{C}{r}+C(1+|x|)\integ{0}{\infty}e^{(\rho-r)s}ds+ \max_{k\in\mathcal{I}}(-C(i,k)).
\end{equation}
Therefore if $r>\rho$ we have
 $$\overline{V}^{\;ij}(x)\leq C(1+|x|).$$

On the other hand, by considering
the particular strategy $\tilde{\delta}=(\tilde{\tau}_m,\tilde{\xi}_m)$ given by $\tilde{\tau}_m=\infty, \tilde{\xi}_m=i$ for all $m\geq 1$, and for $\epsilon>0$ there exists a strategy $\nu\in \mathcal{B}^j$ such that
\begin{equation}
\overline{V}^{\;ij}(x) \geq \mathbb{E}\bigg[\integ{0}{\infty}e^{-rs}f^{a_sb_s}(X_{s}^x)ds+\sum\limits_{n\geq 1}e^{-r\rho_n} \chi(\eta_{n-1},\eta_m)\bigg].
\end{equation}
For all $N\geq1,$ we have
$$\sum_{n=1}^{N}e^{-r\rho_n}\chi(\eta_{n-1},\eta_n)\geq \min_{l\in\mathcal{I}}(\chi(j,l)).$$
Then, we get
\begin{equation}
\overline{V}^{\;ij}(x) \geq \mathbb{E}\bigg[\integ{0}{\infty}e^{-rs}f^{a_sb_s}(X_{s}^x)ds+\min_{l\in\mathcal{I}}(\chi(j,l))\bigg].
\end{equation}
As a consequence, there exists $C$ and $\rho$ such that
\begin{equation}
\overline{V}^{\;ij}(x) \geq -\frac{C}{r}-C(1+|x|)\integ{0}{\infty}e^{(\rho-r)s}ds+\min_{l\in\mathcal{I}}(\chi(j,l)).
\end{equation}
Therefore if $r>\rho$, we have
$$\overline{V}^{\;ij}(x)\geq -C(1+|x|).$$
from which we deduce the thesis.\qquad $\Box$
\begin{lem}
 There exists some positive constant $\rho$ such that for\, $r > \rho$, $x,x'\in\mathbb{R}^m , \delta\in \mathcal{A}^i$ and $\nu\in\mathcal{B}^j$ we
have
$$|J(x,\delta,\nu)-J(x',\delta,\nu)|+|\overline{V}^{\;ij}(x)-\overline{V}^{\;ij}(x')|+|\underline{V}^{ij}(x)-\underline{V}^{ij}(x')|\leq C|x-x'|,$$
for some positive constant C.
\end{lem}
$Proof.$ It is enough to show that the conclusion holds for the gain functional $J$.\\
For every $x,x'\in \mathbb{R}^m, u\in\mathcal{A}^i, \nu\in\mathcal{B}^j$, we have
\begin{multline*}
\qquad\qquad |J(x,\delta,\nu)-J(x',\delta,\nu)|\leq \mathbb{E}\bigg[\int_{0}^{\infty}e^{-rs}|f^{a_sb_s}(X_s^{x})-f^{a_sb_s}(X_s^{x'})|ds\bigg].\qquad\qquad\qquad
\qquad\qquad
\end{multline*}
By a standard estimate for the SDE applying Ito's formula to $|X^x-X^{x'}|^2$ and using Gronwall's lemma, we then obtain from the Lipschitz condition on $b^{ij}, \sigma^{ij}$, the
following inequality
$$E[|X_t^x-X_t^{x'}|\,]\leq e^{\rho t}|x-x'| \qquad\forall x, x'\in\mathbb{R}^m, t\geq 0.$$
From the Lipschitz condition on $f^{ij}$, we deduce
\begin{multline*}
\qquad\qquad |J(x,\delta,\nu)-J(x',\delta,\nu)|\leq C\mathbb{E}\bigg[\int_{0}^{\infty}e^{-rs}|X_s^{x}-X_s^{x'}|ds\bigg]\qquad\qquad\qquad
\qquad\qquad\\ \leq C\mathbb{E}\bigg[\int_{0}^{\infty}e^{-rs}e^{\rho s}|x-x'|ds \bigg]\leq C |x-x'|,\qquad\qquad
\end{multline*}
for $r > \rho$.  This ends the proof.\qquad\qquad $\Box$

In the sequel, we shall assume that $r$ is large enough, which ensures that the value functions  are Lipschitz
continuous. We now present the dynamic programming
principle of  the zero-sum switching game.
\begin{theo}
For any $x\in\mathbb{R}^m$ and  $i,j\in\mathcal{I}$, we have
\begin{equation}
\begin{array}{ll}
\overline{V}^{\;ij}(x)=\inf\limits_{\beta\in {\Delta^j}}\sup\limits_{\delta\in\mathcal{A}^i}\mathbb{E}\bigg[\integ{0}{\theta}e^{-rs}f^{a_sb_s}(X^x_{s})ds-
\sum\limits_{\tau_m\leq \theta}e^{-r\tau_m}C(\xi_{m-1},\xi_m) \\ \qquad\qquad\qquad\qquad\qquad\qquad +\sum\limits_{\rho_n\leq \theta}e^{-r\rho_n}\chi(\eta_{n-1},\eta_n)+e^{-r\theta}\overline{V}^{\;a_\theta b_\theta}(X^x_\theta)\bigg].
\end{array}
\end{equation}
\begin{equation}
\begin{array}{ll}
\underline{V}^{ij}(x)=\sup\limits_{\alpha\in \Gamma^i}\inf\limits_{\nu\in\mathcal{B}^j}\mathbb{E}\bigg[\integ{0}{\theta}e^{-rs}f^{a_sb_s}(X^x_{s})ds-
\sum\limits_{\tau_m\leq \theta}e^{-r\tau_m}C(\xi_{m-1},\xi_m)\qquad \\ \qquad\qquad\qquad\qquad\qquad\qquad+\sum\limits_{\rho_n\leq \theta}e^{-r\rho_n}\chi(\eta_{n-1},\eta_n)+e^{-r\theta}\underline{V}^{a_\theta b_\theta}(X^x_\theta)\bigg],
\end{array}
\end{equation}
where $\theta$ is any stopping time.
\end{theo}
This principle was
formally stated in \cite{[BL]}  and proved rigorously for the finite horizon case in \cite{[TY]}. The arguments
for the infinite horizon case may be adapted in a straightforward way.
\section{Isaacs' system of
quasi-variational inequalities}
\no

In the present section  we prove that
the two value functions are viscosity solutions to (\ref{eq:HJBI1}) and (\ref{eq:HJBI2}). But we make the proof only for (\ref{eq:HJBI2}), the other case is analogous. We begin with the following lemma.
\begin{lem}
\label{lem2}
For any $i,j\in\mathcal{I}$, the lower and
upper value functions satisfy the following properties: for all $x \in\mathbb{R}^m$,
\begin{equation}
\label{Hinf}
V^{ij}(x)\leq N^{ij}[V](x)
\end{equation}
\begin{equation}
\label{Hsup}
M^{ij}[V](x)\leq V^{ij}(x)
\end{equation}
\end{lem}
$Proof:$
We make the proof only for the lower value function $\overline{V}^{\;ij}$, the other case being analogous.
For every $j,\tilde{j}\in\mathcal{I}, j\ne \tilde{j}$ and $\tilde{\nu}\in\mathcal{B}^{\tilde{j}}$ we define $\nu\in\mathcal{B}^j$ by
$$\tilde{\rho}_{l-1}=\rho_l,\; \rho_{0}=0,\; \tilde{\eta}_{l-1}=\eta_l,\; \eta_{0}=j.$$
Note that $\rho_0=\rho_1=0.$ Let $\tilde{\nu}:=\{\tilde{\rho}_l,\tilde{\eta}_l\}$ and $\nu:=\{\rho_l,\eta_l\}$, then
\begin{equation*}
\begin{array}{ll}
 \overline{V}^{\;ij}(x)\leq \sup\limits_{\delta\in\mathcal{A}^i}J(x,\delta,\nu)= \sup\limits_{\delta\in\mathcal{A}^i}\big[J(x,\delta,\tilde{\nu})+\chi(j,\tilde{j})\big],
\end{array}
\end{equation*}
from which we deduce that the following inequality holds:
\begin{equation*}
 \overline{V}^{\;ij}(x)\leq \overline{V}^{\;i\tilde{j}}(x)+\chi(j,\tilde{j}).
\end{equation*}
Therefore
 \begin{equation*}
 \overline{V}^{\;ij}(x)\leq \inf_{\tilde{j}\ne j}\{\overline{V}^{\;i\tilde{j}}(x)+\chi(j,\tilde{j})\}.
\end{equation*}
In a similar way we can prove $(\ref{Hsup})$, hence deducing the thesis. \qquad\qquad$\Box$
\begin{theo}
For any $i,j\in\mathcal{I}$, the lower value function $\overline{V}^{\;ij}$ and the upper value function $\underline{V}^{ij}$ are viscosity solutions of (\ref{eq:HJBI1}).
\end{theo}
$Proof.$
We state the proof only for lower value function, the other case is analogous. First, we prove
the supersolution property. Fix $i,j \in\mathcal{I}$ and let $\bar{x}\in \R^m, \varphi^{ij}\in C^2(\R^m)$ such that $\bar{x}$ is a minimum of $\overline{V}^{\;ij}-\varphi^{ij}$ with $\overline{V}^{\;ij}(\bar{x})=\varphi^{ij}(\bar{x})$.
We have, by Lemma \ref{lem2}
$$\overline{V}^{\;ij}(\bar{x})- N^{ij}[\overline{V}](\bar{x})\leq 0.$$
If
$$\overline{V}^{\;ij}(\bar{x})- N^{ij}[\overline{V}](\bar{x})=0,$$
then we are done. Now suppose
$$\overline{V}^{\;ij}(\bar{x})- N^{ij}[\overline{V}](\bar{x})<-2\epsilon
<0,$$
we prove by contradiction that
$$r\varphi^{ij}(\bar{x})-\mathcal{A}\varphi^{ij}(\bar{x})-f^{ij}(\bar{x})\geq 0.$$
Suppose otherwise, i.e., $r\varphi^{ij}(\bar{x})-\mathcal{A}\varphi^{ij}(\bar{x})-f^{ij}(\bar{x})< 0$. Then without loss of
generality we can assume that $r\varphi^{ij}(\bar{x})-\mathcal{A}\varphi^{ij}(\bar{x})-f^{ij}(\bar{x})< 0$ and $\overline{V}^{\;ij}(x)- N^{ij}[\overline{V}](x)<-\epsilon
<0$ on $B(\bar{x},\delta)$. \\
Define the stopping time $\tau$ by
$$\tau=\inf\{t\geq 0: X_{t}\notin B(\bar{x},\delta)\}.$$
Let $\epsilon_1>0$, using the dynamic programming principle, we deduce the existence of a strategy $\beta\in\Delta^j$ such that
\begin{equation*}
\begin{array}{llllll}
\overline{V}^{\;ij}(\bar{x})\geq \mathbb{E}\bigg[\integ{0}{\tau\wedge\rho_1}e^{-rs}f^{ij}(X_{s}^{\bar{x}})ds+\ind_{[\rho_1\leq \tau]}e^{-r\rho_1}\big(\chi(j,\eta_1)+\overline{V}^{\;i\eta_1}(X_{\rho_1}^{\bar{x}})\big)\bigg]\\ \\ \qquad\qquad+\mathbb{E}\bigg[\ind_{[\rho_1 > \tau]}e^{-r\tau}\overline{V}^{\;ij}(X_{\tau}^{\bar{x}})\bigg]-\epsilon_1
\\ \\ \qquad\quad\; \geq \mathbb{E}\bigg[\integ{0}{\tau\wedge\rho_1}e^{-rs}f^{ij}(X_{s}^{\bar{x}})ds+\ind_{[\rho_1\leq \tau]} e^{-r\rho_1}N^{ij}[\overline{V}](X_{\rho_1}^{\bar{x}})\bigg]+\mathbb{E}\bigg[\ind_{[\rho_1 > \tau]}e^{-r\tau}\overline{V}^{\;ij}(X_{\tau}^{\bar{x}})\bigg]-\epsilon_1
\\ \\ \qquad\quad\; \geq \mathbb{E}\bigg[\integ{0}{\tau\wedge\rho_1}e^{-rs}f^{ij}(X_{s}^{\bar{x}})ds+e^{-r(\rho_1\wedge\tau)}\overline{V}^{\;ij}(X_{\rho_1\wedge\tau}^{\bar{x}})\bigg]-\epsilon_1.
\end{array}
\end{equation*}
Therefore, without loss of generality, we only need to consider $\beta\in\Delta^j$ such that $\rho_1 > \tau$. Then
\begin{equation*}
\begin{array}{llllll}
\varphi^{ij}(\bar{x})=\overline{V}^{\;ij}(\bar{x})\geq \mathbb{E}\bigg[\integ{0}{\tau}e^{-rs}f^{ij}(X_{s}^{\bar{x}})ds+ e^{-r\tau}\overline{V}^{\;ij}(X_{\tau}^{\bar{x}})\bigg]-\epsilon_1\qquad\qquad\qquad
\\  \\ \qquad\qquad\quad\quad\quad\geq \mathbb{E}\bigg[\integ{0}{\tau}e^{-rs}f^{ij}(X_{s}^{\bar{x}})ds+ e^{-r\tau}\varphi^{ij}(X_{\tau}^{\bar{x}})\bigg]-\epsilon_1.
\end{array}
\end{equation*}
By applying It\^o's formula to $e^{-rt}\varphi^{ij}(X^{\bar{x}}_t)$ and plugging into
the last inequality, we obtain
$$
 0\leq \mathbb{E}\bigg[\integ{0}{\tau}e^{-rs}(r\varphi^{ij}-\mathcal{A}\varphi^{ij}-f^{ij})(X_{s}^{\bar{x}})ds\bigg]-\epsilon_1,
$$
which is a contradiction. Therefore, we must have
$
 0\leq (r\varphi^{ij}-\mathcal{A}\varphi^{ij}-f^{ij})(\bar{x})
$.
\\
Thanks to Lemma $\ref{lem2}$, it is enough to show that given
$\bar{x}\in\mathbb{R}^{m}$ such that
$$M^{ij}[\overline{V}](\bar{x})\leq \overline{V}^{\;ij}(\bar{x}).$$
Then
\begin{equation}
\begin{array}{c}
max\Big\{ min\Big[r\varphi^{ij}(\bar{x})-\mathcal{A}\varphi^{ij}(\bar{x})-f^{ij}(\bar{x});\qquad\qquad\qquad\qquad\\ \qquad\qquad\overline{V}^{\;ij}(\bar{x})-M^{ij}[\overline{V}](\bar{x})\Big];\overline{V}^{\;ij}(\bar{x})-N^{ij}[\overline{V}](\bar{x})\Big\}\geq 0.
\end{array}
\end{equation}
Which is the supersolution property. The subsolution property is proved analogously.\qquad$\Box$
\section{Uniqueness of the solution of Isaacs' system of
quasi-variational inequalities}
\no

We prove that the (2.9) admits a unique viscosity solution. We need to make an additional assumption on the switching costs.\\
\textbf{[H5]} (The no free loop property)

For  any sequence of pairs $\{i_p,j_p\}_{p=1}^{N+1}\subset\mathcal{I}\times\mathcal{I}$ such that $(i_{N+1},j_{N+1})=(i_1,j_1),$ and for $ 1\leq p\leq N ,$ either $i_{p+1}=i_p$ or $j_{p+1}=j_p$, we have
$$\sum_{p=1}^{N}\chi(j_p,j_{p+1})-\sum_{p=1}^{N}C(i_p,i_{p+1})\ne 0.$$

We begin with the technical lemma.
\begin{lem}\label{lem4}
Let $\vec{U}:=(U^{ij})_{(i,j)\in\mathcal{I}\times\mathcal{I}}$ and $\vec{V}:=(V^{ij})_{(i,j)\in\mathcal{I}\times\mathcal{I}}$ be a viscosity subsolution and a viscosity supersolution to the  equation (\ref{eq:HJBI1}), respectively. Let $x\in\mathbb{R}^m$ be fixed and $(\bar{i},\bar{j})\in\mathcal{I}\times\mathcal{I}$ such that
\begin{equation}
V^{\bar{i}\bar{j}}(x)-U^{\bar{i}\bar{j}}(x)=\max_{i,j\in\mathcal{I}}\{V^{ij}(x)-U^{ij}(x)\}
\end{equation}
and
\begin{equation}\label{eq111unic}
U^{\bar{i}\bar{j}}(x)\leq M^{\bar{i}\bar{j}}[U](x)\quad \text{or}\quad V^{\bar{i}\bar{j}}(x) \geq N^{\bar{i}\bar{j}}[V](x).
\end{equation}
Then  there exists $(i_0,j_0)\in\mathcal{I}\times\mathcal{I}$ such that
\begin{equation}\label{eq1unic}
V^{\bar{i}\bar{j}}(x)-U^{\bar{i}\bar{j}}(x)= V^{i_0j_0}(x)-U^{i_0j_0}(x)
\end{equation}
and
\begin{equation}\label{eq12unic}
U^{i_0j_0}(x)> M^{i_0j_0}[U](x)\quad \text{and}\quad V^{i_0j_0}(x) < N^{i_0j_0}[V](x).
\end{equation}
\end{lem}
$Proof$.
Suppose that $U^{\bar{i}\bar{j}}(x)\leq M^{\bar{i}\bar{j}}[U](x_{\epsilon})$, then there exists $k\ne\bar{i}$, such that
\begin{equation}\label{UN1}
U^{\bar{i}\bar{j}}(x)\leq U^{k\bar{j}}(x)-C(\bar{i},k).
\end{equation}
Since $\vec{V}$ is a subsolution, it satisfies
\begin{equation}\label{UN2}
V^{\bar{i}\bar{j}}(x)\geq V^{k\bar{j}}(x)-C(\bar{i},k)
\end{equation}
which implies that
$$U^{\bar{i}\bar{j}}(x)-V^{\bar{i}\bar{j}}(x)\leq U^{k\bar{j}}(x)-V^{k\bar{j}}(x).$$
Hence
\begin{equation}\label{UN3}
U^{\bar{i}\bar{j}}(x)-V^{\bar{i}\bar{j}}(x)= U^{k\bar{j}}(x)-V^{k\bar{j}}(x).
\end{equation}
Then (\ref{UN1}), (\ref{UN2}) and (\ref{UN3})  imply
$$U^{\bar{i}\bar{j}}(x)-U^{k\bar{j}}(x)= V^{\bar{i}\bar{j}}(x)-V^{k\bar{j}}(x)=-C(\bar{i},k).$$

Similarly if $V^{\bar{i}\bar{j}}(x) \geq N^{\bar{i}\bar{j}}[V](x)$ hold, then there exists  $l\ne\bar{j}$ such that
\begin{equation*}
U^{\bar{i}\bar{j}}(x)-V^{\bar{i}\bar{j}}(x)= U^{\bar{i}l}(x)-V^{\bar{i}l}(x)
\end{equation*}
and
$$U^{\bar{i}\bar{j}}(x)-U^{\bar{i}l}(x)= V^{\bar{i}\bar{j}}(x)-V^{\bar{i}l}(x)=\chi(\bar{j},l).$$
Now if the new index ($k$ or $l$)  verify (\ref{eq111unic}) , we repeat this reasoning . If this case continues to occur, finally we find a loop $(i_1,j_1),...,(i_{N},j_{N}),(i_{N+1},j_{N+1})=(i_1,j_1)$ such that
$$\sum_{q=1}^{N}\chi(j_q,j_{q+1})-\sum_{q=1}^{N}C(i_q,i_{q+1})=0.$$
Hence we obtain a contradiction to the assumption (\textbf{H5}). Thus (\ref{eq1unic}) and (\ref{eq12unic}) holds for some $(i_0,j_0)$ \quad$\Box$
\begin{theo}
\label{thoUnicitÃ©}
Let $\vec{U}=(U^{ij})_{(i,j)\in\mathcal{I}\times\mathcal{I}}$ (resp. $\vec{V}=(V^{ij})_{(i,j)\in\mathcal{I}\times\mathcal{I}}),$ a family of continuous viscosity subsolutions (resp.
 supersolutions) to (\ref{eq:HJBI1}), and satisfying a linear growth condition. Then, $U^{ij} \leq V^{ij}$ for all $i,j\in\mathcal{I}.$
\end{theo}
$Proof.$
Let us proceed by contradiction. For some $R > 0$ suppose there exists there exists $(\bar{x},\bar{i},\bar{j})\in B_R\times\mathcal{I}\times\mathcal{I} \;(B_R:=\{x\in\mathbb{R}^m; |x|<R \})$ such that
\begin{equation}
\max\limits_{x\in B_R}\max_{i,j}(U^{ij}-V^{ij})(x)=(U^{\bar{i}\bar{j}}-V^{\bar{i}\bar{j}})(\bar{x})=\eta > 0.
\end{equation}
We divide the proof into two steps.

\textbf{Step 1.} Using Lemma \ref{lem4} we derive the existence of $i_0$ and $j_0$ such that
\begin{equation}
(U^{i_0j_0}-V^{i_0j_0})(\bar{x})=\eta > 0.
\end{equation}
and
\begin{equation}
V^{i_0j_0}(x) < N^{i_0j_0}[V](x)\quad \text{and}\quad U^{i_0j_0}(x)> M^{i_0j_0}[U](x)
\end{equation}
For a small $\epsilon>0$, let $\Phi_{\epsilon}$ be the function defined as follows.
\begin{equation}
\begin{array}{ll}
\label{phi}
\Phi_{\epsilon}(x,y)=U^{i_0j_0}(x)-V^{i_0j_0}(y)-\displaystyle\frac{1}{2\epsilon}|x-y|^{2}\qquad\qquad (x,y)\in\mathbb{R}^m\times\mathbb{R}^m.
\end{array}
\end{equation}
Let $(x_\epsilon,y_\epsilon)$  be such that
\begin{equation}
\begin{array}{ll}
\label{phi}
\Phi^{i_\epsilon j_\epsilon}_{\epsilon}(x_\epsilon,y_\epsilon)=\max\limits_{x,y\in B_R}\max\limits_{i,j\in\mathcal{I}}\Phi^{ij}_{\epsilon}(x,y).
\end{array}
\end{equation}
From $2\Phi_{\epsilon}(x_\epsilon,y_\epsilon)\geq\Phi_{\epsilon}(x_\epsilon,x_\epsilon)+\Phi(y_\epsilon,y_\epsilon)$ we have
\begin{equation}
\displaystyle\frac{1}{\epsilon}|x_\epsilon -y_\epsilon|^{2} \leq
(U^{i_0j_0}(x_\epsilon)-U^{i_0 j_0}(y_\epsilon))+(V^{i_0 j_0}(x_\epsilon)-V^{i_0 j_0}(y_\epsilon)),
\end{equation}
and consequently $\displaystyle\frac{1}{\epsilon}|x_\epsilon -y_\epsilon|^{2}$ is bounded,
and as $\epsilon\rightarrow 0$, $|x_\epsilon -y_\epsilon|\rightarrow 0$. Since
$U^{i_\epsilon j_\epsilon}$ and $V^{i_\epsilon j_\epsilon}$ are uniformly continuous, then $\displaystyle\frac{1}{2\epsilon}|x_\epsilon
-y_\epsilon|^{2}\rightarrow 0$ as
$\epsilon\rightarrow 0.$\\
Next, recalling that $U^{i_0 j_0}$ and $V^{i_0 j_0}$ are continuous, then, for $\epsilon$ small enough and at least for a subsequence which we still index by $\epsilon$, we obtain
\begin{equation}
V^{i_0j_0}(y_{\epsilon}) < N^{i_0j_0}[V](y_{\epsilon})\quad \text{and}\quad U^{i_0j_0}(x_{\epsilon})> M^{i_0j_0}[U](x_{\epsilon})
\end{equation}
\textbf{Step 2.}
Let us denote
\begin{equation}
\varphi_{\epsilon}(x,y)=\displaystyle\frac{1}{2\epsilon}|x-y|^{2}.
\end{equation}
Then we have: \be \left\{
\begin{array}{lllllll}\label{derive}
D_x\varphi_{\epsilon}(x,y)= \displaystyle\frac{1}{\epsilon}(x-y), \\
D_y\varphi_{\epsilon}(x,y)= -\displaystyle\frac{1}{\epsilon}(x-y)\\
\\
B(x,y)=D^2_{x,y}\varphi_\epsilon(x,y)=\displaystyle\frac{1}{\epsilon}\begin{pmatrix}
   I & -I \\
   -I & I
\end{pmatrix}.
\end{array}
\right. \ee
Then applying the result by Crandall et al. (Theorem 3.2, \cite{[CIL]}) to the function
$$U^{i_0 j_0}(x)-V^{i_0 j_0 }(y)-\varphi_\epsilon(x,y)$$
at the point $(x_\epsilon,y_\epsilon)$, for any $\epsilon_1>0$, we can find  $X,Y\in \mathbb{S}_m,$ such that:
 \be \left\{
\begin{array}{lllllll}\label{derive1}
\big(\displaystyle\frac{1}{\epsilon}(x_\epsilon-y_\epsilon),X\big)\in J^{2,+}(U^{i_0  j_0}(x_\epsilon)), \\ \\
\big(\displaystyle\frac{1}{\epsilon}(x_\epsilon-y_\epsilon),Y\big)\in J^{2,-}(V^{i_0j_0 }(y_\epsilon)),\\
\\
-\big(\displaystyle\frac{1}{\epsilon_1}+\|B(x_\epsilon,y_\epsilon)\|\big)I\leq \begin{pmatrix}
   X & 0 \\
   0 & -Y
\end{pmatrix}\leq B(x_\epsilon,y_\epsilon)+\epsilon_1B(x_\epsilon,y_\epsilon)^2.
\end{array}
\right. \ee
 Then by definition of
viscosity solution, we get:
\begin{equation}\begin{array}{lll}\label{vis_sub1}
rU^{i_0j_0}(x_\epsilon)
-\langle\displaystyle\frac{1}{\epsilon}(x_\epsilon-y_\epsilon),\\ \qquad\qquad
b(x_\epsilon)\rangle-\displaystyle\frac{1}{2}tr[\sigma^*(x_\epsilon)X\sigma(x_\epsilon)]-f^{i_0 j_0 }(x_\epsilon)]\leq0,
\end{array}\end{equation}
 and
\begin{equation}\begin{array}{l}\label{vis_sub11}
rV^{i_0 j_0 }(y_\epsilon)
-\langle\displaystyle\frac{1}{\epsilon}(x_\epsilon-y_\epsilon),\\ \qquad\qquad\qquad\qquad
b(y_\epsilon)\rangle-\displaystyle\frac{1}{2}tr[\sigma^*(y_\epsilon)Y\sigma(y_\epsilon)]-f^{i_0 j_0 }(y_\epsilon)\geq0,
\end{array}\end{equation}
which implies that:
\begin{equation}
\begin{array}{llllll}
\label{viscder} &rU^{i_0 j_0 }(x_\epsilon)-rV^{i_0 j_0 }(y_\epsilon)\\ \\ & \leq
\langle\displaystyle\frac{1}{\epsilon}(x_\epsilon-y_\epsilon) ,
b(x_\epsilon)-b(y_\epsilon)\rangle
+\displaystyle\frac{1}{2}tr[\sigma^*(x_\epsilon)X\sigma(x_\epsilon)-\sigma^*(y_\epsilon)Y\sigma(y_\epsilon)]\\ \\
&
+f^{i_0 j_0}(x_\epsilon)-f^{i_0 j_0 }(y_\epsilon).
\end{array}
\end{equation}
we have :
\begin{equation}
B+\epsilon_1B^2\leq \frac{\epsilon+\epsilon_1}{\epsilon^2}\begin{pmatrix}
   I & -I \\
   -I & I
\end{pmatrix},
\end{equation}
where $C$ which hereafter may change from line to line. Choosing now $\epsilon_1=\epsilon,$ yields the relation
\begin{equation}
\label{equaB}
B+\epsilon_1B^2\leq \frac{2}{\epsilon}\begin{pmatrix}
   I & -I \\
   -I & I
\end{pmatrix}.
\end{equation}
Now, from (\textbf{H1}), (\ref{derive1}) and (\ref{equaB}) we get:
$$\displaystyle\frac{1}{2}tr[\sigma^*(x_\epsilon)X\sigma(x_\epsilon)-\sigma^*(y_\epsilon)Y\sigma(y_\epsilon)]\leq \frac{C}{\epsilon}|x_\epsilon-y_\epsilon|^2.$$
Next
$$
\langle\frac{1}{\epsilon}(x_\epsilon-y_\epsilon),b(x_\epsilon)-b(y_\epsilon)\rangle
\leq \frac{C^2}{\epsilon}|x_\epsilon - y_\epsilon|^{2}.$$
So that by plugging into (\ref{viscder}) we obtain:
\begin{equation}
\begin{array}{llllll}
\label{viscder11}
rU^{i_0j_0}(x_\epsilon)-rV^{i_0 j_0}(y_\epsilon)\\ \\
\qquad\qquad\leq \displaystyle\frac{C}{\epsilon}|x_\epsilon-y_\epsilon|^2+\displaystyle\frac{C^2}{\epsilon}|x_\epsilon - y_\epsilon|^{2}+f^{i_0 j_0}(x_\epsilon)-f^{i_0 j_0}(y_\epsilon).
\end{array}
\end{equation}By sending $\epsilon \rightarrow0$, and taking into account of the continuity of $f^{i_0 j_0}$, we obtain $\eta \leq 0$ which is a contradiction. Thus, $U^{ij}\leq V^{ij}$, for any $(i,j)\in\mathcal{I}\times\mathcal{I}$, which is the desired result. \qquad$\Box$
\begin{cor}
The lower and
upper value functions coincide, and the value function of the stochastic differential
game is given by $V^{ij}(x):=\overline{V}^{\;ij}(x)=\underline{V}^{ij}(x)$ for every $i,j\in\mathcal{I}$ and $x\in\mathbb{R}^m.$ As a consequence the two equations (\ref{eq:HJBI1}) and (\ref{eq:HJBI2}) coincide.
\end{cor}
\section{A verification theorem}
\no

In this section, we present a verification theorem which gives an optimal
strategies of our zero-sum stochastic differential game.

We suppose that a classical solution of (\ref{eq:HJBI1}) exists, denoted by $(V^{ij})_{(i,j)\in\mathcal{I}\times\mathcal{I}}$. Then for each $i,j\in\mathcal{I}, V^{ij}$ separates the space $\mathbb{R}^m$ into
four regions:
\begin{equation*}
\begin{array}{ll}
\mathcal{C}:=\{x\in\mathbb{R}^n: V^{ij}(x)<N^{ij}[V](x) ; V^{ij}(x)> M^{ij}[V](x)\; \mbox{and}\; rV^{ij}(x)-\mathcal{A}V^{ij}(x)-f^{ij}(x)=0\}\\ \\
\mathcal{I}_1:=\{x\in\mathbb{R}^n: V^{ij}(x)<N^{ij}[V](x) ; V^{ij}(x)= M^{ij}[V](x) \;\mbox{and}\;  rV^{ij}(x)-\mathcal{A}V^{ij}(x)-f^{ij}(x)\geq 0\}\\ \\
\mathcal{I}_2:=\{x\in\mathbb{R}^n: V^{ij}(x)=N^{ij}[V](x) ; V^{ij}(x)> M^{ij}[V](x) \;\mbox{and}\; rV^{ij}(x)-\mathcal{A}V^{ij}(x)-f^{ij}(x)\leq 0\}\\ \\
\mathcal{I}_3:=\{x\in\mathbb{R}^n: V^{ij}(x)=N^{ij}[V](x) \quad\mbox{and}\quad  V^{ij}(x)= M^{ij}[V](x)\}
\end{array}
\end{equation*}

Let us define the strategies $\delta^*:=(\tau^*_m,\xi^*_m)_{m\geq0}$ (resp. $\nu^* := (\rho^*_n,\eta^*_n)_{n\geq0})$ as follows:
$$\tau^*_0=0,\xi^*_0=i \;(\text{resp}.\; \rho^*_0=0,\eta^*_0=j)$$
and for any $m \geq 1$,
 \begin{equation*}
\left\{
\begin{array}{ll}
 \tau^*_m=\left\{\begin{array}{ll}\text{inf}\{s\geq\tau^*_{m-1}, V^{\xi^*_{m-1}b_s}(X_s)=\text{max}_{k\ne \xi^*_{m-1}}\{V^{k\, b_s}(X_s)-C(\xi^*_{m-1},k)\}\}\\ \\  +\infty \quad\text{if the above set is empty}\end{array}\right.
 \\ \\ \text{and} \\ \\
\xi^*_m=\left\{\begin{array}{ll}\text{max}\{k\ne \xi^*_{m-1},V^{k\,a_{\tau^*_m}}(X_{\tau^*_m})-C(\xi^*_{m-1},k)\}\quad\mbox{if}\quad \tau^*_m<+\infty\\ \\  \xi^*_{m-1} \quad\mbox{if}\quad \tau^*_m=+\infty \end{array}\right.
\end{array}
\right. \end{equation*}
\bigg(resp.
 \begin{equation*} \left\{
\begin{array}{ll}
 \rho^*_m=\left\{\begin{array}{ll}
 \inf\{s\geq\rho^*_{m-1}, V^{a_s\eta^*_{m-1}}(X_s)=\text{min}_{l\ne \eta^*_{m-1}}\{V^{a_sl}(X_s)+\chi(\eta^*_{m-1},l)\}\} \\ \\ +\infty \quad\text{if the above set is empty}\end{array}\right.
 \\ \\ \text{and} \\ \\
\eta^*_m=\left\{
\begin{array}{ll}
\min\{l\ne\eta^*_{m-1},V^{a_{\rho^*_m}l}(X_{\rho^*_m})+\chi(\eta^*_{m-1},l)\}\quad\text{if}\quad \rho^*_m < +\infty \\ \\ \eta^*_{m-1} \quad\text{if}\quad \rho^*_m=+\infty
\end{array}\right.
\bigg)\end{array}
\right. \end{equation*}

We are now ready to present the verification theorem for our switching game.
\begin{theo}
For each $i,j\in\mathcal{I}$, $x\in\mathbb{R}^m$, and assume that $(\delta^*,\nu^*)\in\mathcal{A}^{i}\times\mathcal{B}^{j}$. Then we have
 $V^{ij}(x)=J(x,\delta^*,\nu^*).$
\end{theo}
$Proof:$
First for each $\delta:=(\tau_m,\xi_m)_{m\geq0}\in\mathcal{A}^{i}$ and $\nu := (\rho_n,\eta_n)_{n\geq0}\in\mathcal{B}^{j}$ we define $\theta=(\theta_k)$ an increasing sequence of stopping times by
\begin{equation*}
\begin{array}{ll}
\qquad\qquad\theta_0=\tau_0=\rho_0=0\\ \\
\qquad\qquad\theta_1=\min(\tau_1,\rho_1)\\ \\
\qquad\qquad\theta_2=\min(\tau_1 \ind_{[\tau_1>\theta_{1}]}+\tau_2\ind_{[\tau_1\leq\theta_{1}]},\rho_1\ind_{[\rho_1>\theta_{1}]}+\rho_2\ind_{[\rho_1\leq\theta_{1}]}),\\ \\
\qquad\qquad\theta_3=\min(\tau_1\ind_{[\tau_1>\theta_{2}]}+\tau_3\ind_{[\tau_1\leq\theta_{2}]},\rho_1\ind_{[\rho_1>\theta_{2}]}+\rho_3\ind_{[\rho_1\leq\theta_{2}]}\\
\qquad\qquad\qquad\qquad\quad,\tau_2\ind_{[\tau_2>\theta_{2}]}+\tau_3\ind_{[\tau_2\leq\theta_{2}]},\rho_2\ind_{[\rho_2>\theta_{2}]}+\rho_3\ind_{[\rho_2\leq\theta_{2}]}),\\ \\
\qquad\qquad\qquad\qquad\qquad.....................................................\\ \\
\qquad\qquad\theta_k=\min(\tau_1\ind_{[\tau_1>\theta_{k-1}]}+\tau_k\ind_{[\tau_1\leq\theta_{k-1}]},\rho_1\ind_{[\rho_1>\theta_{k-1}]}+\rho_k\ind_{[\rho_1\leq\theta_{k-1}]},...

\\ \qquad\qquad\qquad\qquad,\tau_{k-1}\ind_{[\tau_{k-1}>\theta_{k-1}]}+\tau_k\ind_{[\tau_{k-1}\leq\theta_{k-1}]} ,\rho_{k-1}\ind_{[\rho_{k-1}>\theta_{k-1}]}+\rho_k\ind_{[\rho_{k-1}\leq\theta_{k-1}]}).
\end{array}
\end{equation*}
Then the cost functional is rewritten as:
\begin{equation}
\label{eqver0}
\begin{array}{ll}
J(x,\delta,\nu)=\displaystyle\sum\limits_{k\geq 1}\mathbb{E}\bigg[\integ{\theta_{k-1}}{\theta_{k}}e^{-rs}f^{a_sb_s}(X_{s})ds-
\sum\limits_{m\geq 1}e^{-r\tau_m} C(\xi_{m-1},\xi_m)\ind_{[\tau_m=\theta_{k}]} \\ \qquad\qquad\qquad\qquad +\displaystyle\sum\limits_{n\geq 1}e^{-r\rho_n}\chi(\eta_{n-1},\eta_n)\ind_{[\rho_n=\theta_{k}]}\bigg].
\end{array}
\end{equation}
Now, let $(\theta^*)_k$ associated with $\delta^*$ and $\nu^*$, then when $\theta^*_{k-1} < t < \theta^*_{k}$ we have
 \begin{equation}
rV^{a_sb_t}(X_t)-\mathcal{A}V^{a_tb_t}(X_t)=f^{a_tb_t}
 \end{equation}
Then by It\^o's formula (see, e.g.,
Sect. IV.45 of \cite{[RW]}), we obtain
\begin{equation}
\begin{array}{ll}
\label{eqver00}
\mathbb{E}\bigg[\integ{\theta^*_{k-1}}{\theta^*_{k}}e^{-rs}f^{a_sb_s}(X_{s})ds\bigg]=\mathbb{E}\bigg[\integ{\theta^*_{k-1}}{\theta^*_{k}}rV^{a_sb_s}(X_s)-\mathcal{A}V^{a_sb_s}(X_s)ds\bigg]\\ \\ \qquad\qquad\qquad\qquad\qquad\quad\; =\mathbb{E}\big[e^{-r\theta^*_{k-1}}V^{a_{\theta^*_{k-1}}b_{\theta^*_{k-1}}}(X_{\theta^*_{k-1}})-e^{-r\theta^*_{k}}V^{a_{\theta^*_{k-1}}b_{\theta^*_{k-1}}}(X_{\theta^*_{k}})\big]
\end{array}
\end{equation}
Substituting this into (\ref{eqver0}), we obtain
\begin{equation}
\begin{array}{ll}
\label{eqver1}
J(x,\delta^*,\nu^*)=\mathbb{E}\bigg[V^{ij}(x)+\displaystyle\sum_{k\geq1}\big[\{-V^{a_{\theta^*_{k-1}}b_{\theta^*_{k-1}}}(X_{\theta^*_{k}})+V^{a_{\theta^*_{k}}b_{\theta^*_{k}}}(X_{\theta^*_{k}}) \\ \qquad\qquad-
\displaystyle\sum\limits_{m\geq 1}C(\xi^*_{m-1},\xi^*_m)\ind_{[\tau^*_m=\theta^*_{k}]} +\displaystyle\sum\limits_{n\geq 1}\chi(\eta^*_{n-1},\eta^*_n)\ind_{[\rho^*_n=\theta^*_{k}]}\}e^{-r\theta^*_{k}}\big]\bigg].
\end{array}
\end{equation}
We now estimate the term in the right-hand side of (\ref{eqver1}) for each $k \geq 1$.\\
\begin{itemize}
\item
  If $\theta^*_{k}= \tau^*_m$ for some $m\geq 1$ we have
\begin{equation}
\label{esti1}
V^{\xi^*_{m}b_{\theta^*_{k-1}}}(X_{\theta^*_{k}})-V^{\xi^*_{m-1}b_{\theta^*_{k-1}}}(X_{\theta^*_{k}})-C(\xi^*_{m-1},\xi^*_m)=0
\end{equation}
\item
  If $\theta^*_{k}= \rho^*_n$ for some $n\geq 1$ we have
\begin{equation}
\label{esti2}
V^{a_{\theta^*_{k-1}}\eta^*_{n}}(X_{\theta^*_{k}})-V^{a_{\theta^*_{k-1}}\eta^*_{n-1}}(X_{\theta^*_{k}})+\chi(\eta^*_{n-1},\eta^*_n)=0
\end{equation}
\item
  If $\theta^*_{k}=\tau^*_m=\rho^*_n$ for some  $m,n\geq 1$ we have
\begin{equation}
\label{esti3}
V^{\xi^*_{m}\eta^*_{n}}(X_{\theta^*_{k}})-V^{\xi^*_{m-1}\eta^*_{n-1}}(X_{\theta^*_{k}})+\chi(\eta^*_{n-1},\eta^*_n)-
C(\xi^*_{m-1},\xi^*_m)=0.
\end{equation}
\end{itemize}
By (\ref{eqver1}) and the estimates in (\ref{esti1})-(\ref{esti3}) it follows that
$$J(x,\delta^*,\nu^*)=V^{ij}(x)\qquad\qquad\Box$$
\section{Numerical results}
\no

In this section, we present the results of some numerical simulations on the game by
using MATLAB, here we apply the policy iteration algorithm for solving numerically(\ref{eq:HJBI1}).\\
In particular, for the numerical example we consider a two-regime switching problem where the diffusion is independent of the regime and follows a geometric Brownian motion, i.e $b^{ij}(x)=bx$ and $\sigma^{ij}(x)=\sigma x$  for
some $b\in\mathbb{R}$ and $\sigma > 0$ .\\ We consider the following game problem:
\begin{equation}
\begin{array}{ll}
\qquad\quad\mathcal{I}=\{1,2\}\\ \\
f^{11}=5x,\qquad f^{12}=x,\\ \\
f^{21}=-x,\qquad f^{22}=-4x,\\ \\
r=0.15,\quad \sigma=0.2, \quad b=0.01,
\end{array}
\end{equation}

In the flowing figures, we plot the value functions for different switching costs.
\begin{figure}[ht]\label{fig1}
  \begin{subfigure}[b]{0.5\linewidth}
    \centering
    \includegraphics[width=0.9\linewidth]{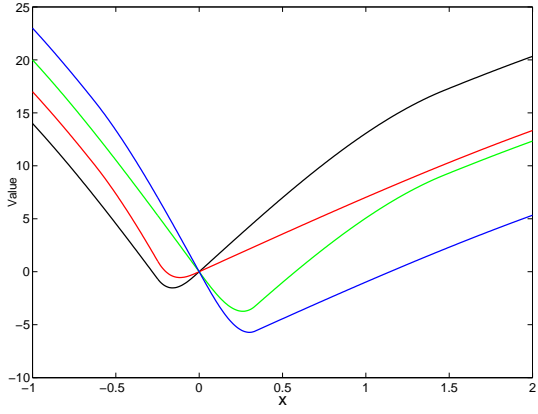}
    \caption{\scriptsize{$C(1,2)=2,C(2,1)=5,\chi(1,2)=2,\chi(2,1)=5$}}
    \label{fig7:a}
    \vspace{4ex}
  \end{subfigure}
  \begin{subfigure}[b]{0.5\linewidth}
    \centering
    \includegraphics[width=0.9\linewidth]{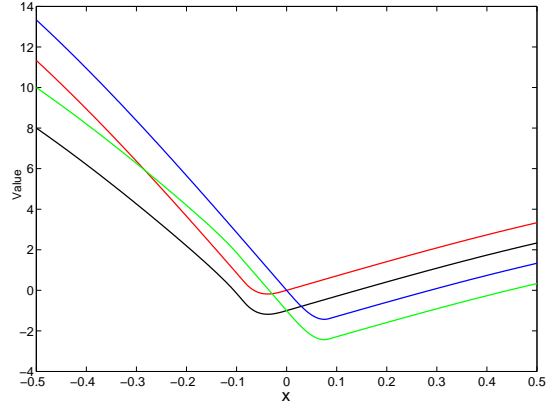}
    \caption{\scriptsize{$C(1,2)=5,C(2,1)=-1,\chi(1,2)=2,\chi(2,1)=2$}}
    \label{fig7:b}
    \vspace{4ex}
  \end{subfigure}
  \begin{subfigure}[b]{0.5\linewidth}
    \centering
    \includegraphics[width=0.9\linewidth]{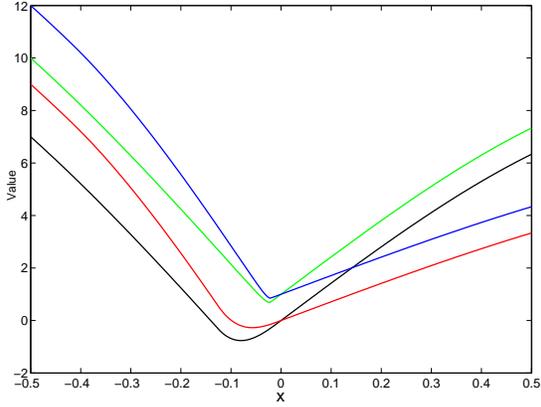}
    \caption{\scriptsize{$C(1,2)=2,C(2,1)=2,\chi(1,2)=-1,\chi(2,1)=5$}}
    \label{fig7:c}
  \end{subfigure}
  \begin{subfigure}[b]{0.5\linewidth}
    \centering
    \includegraphics[width=0.9\linewidth]{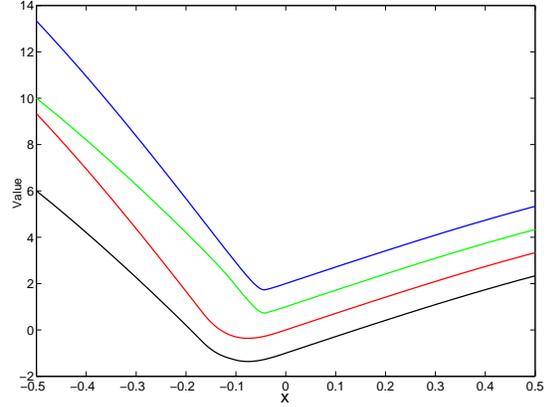}
    \caption{\scriptsize{$C(1,2)=4,C(2,1)=-2,\chi(1,2)=-1,\chi(2,1)=3$}}
    \label{fig7:d}
  \end{subfigure}
  \caption{Value functions:$V^{11}$(black),$V^{12}$(red),$V^{21}$(green),$V^{22}$(blue)}.
  \label{fig7}
\end{figure}
\newpage
Next, we plot the optimal strategies of our zero-sum stochastic differential game. We focus on exemple $(a)$ in figure $1$ . First let  $\mathcal{I}_{ij\to kl}=\{x\in\mathbb{R}: V^{ij}(x)=V^{kl}(x)-C(i,k)+\chi(j,l)\}$. Then by figure $1$ (exemple $(a))$ we get
\begin{equation}
\begin{array}{ll}
\mathcal{I}_{11\to 21}=(-\infty,-0.25] \qquad\quad \mathcal{I}_{11\to 12}=[1.46,+\infty)\\\mathcal{I}_{12\to 21}=(-\infty,-0.62] \qquad\quad \mathcal{I}_{12\to 22}=[-0.62,-0.25]\\\mathcal{I}_{21\to 11}=[0.33,1.46] \qquad\quad\;\;\; \mathcal{I}_{21\to 12}=[1.46,+\infty)\\\mathcal{I}_{22\to 21}=(-\infty,-0.62] \qquad\quad \mathcal{I}_{22\to 12}=[0.33,+\infty)
\end{array}
\end{equation}
\begin{figure}[t]
{\includegraphics[width=1.02\textwidth]{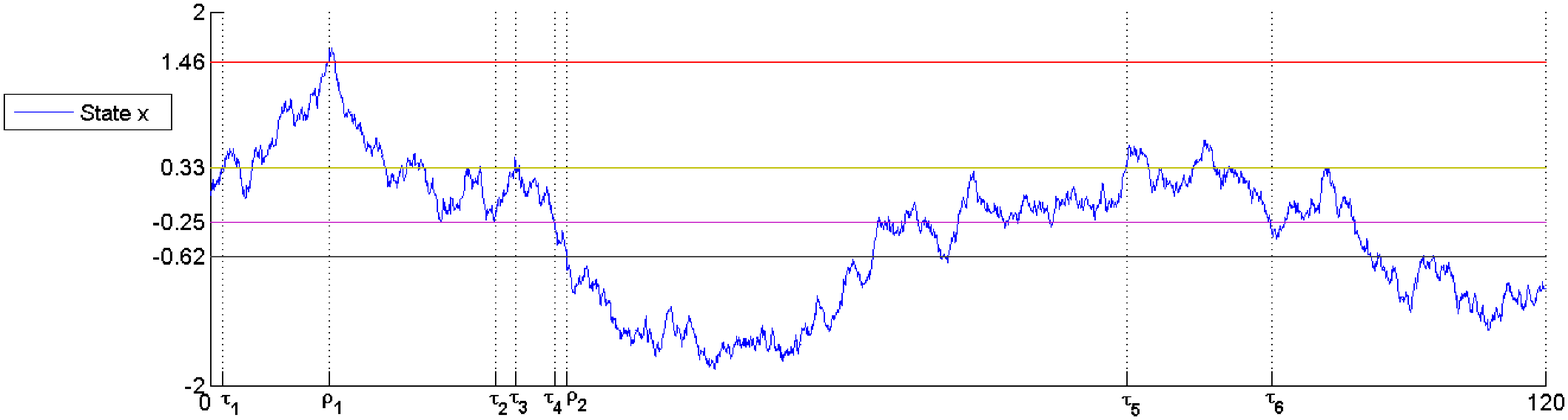}}
{\includegraphics[width=1.02\textwidth]{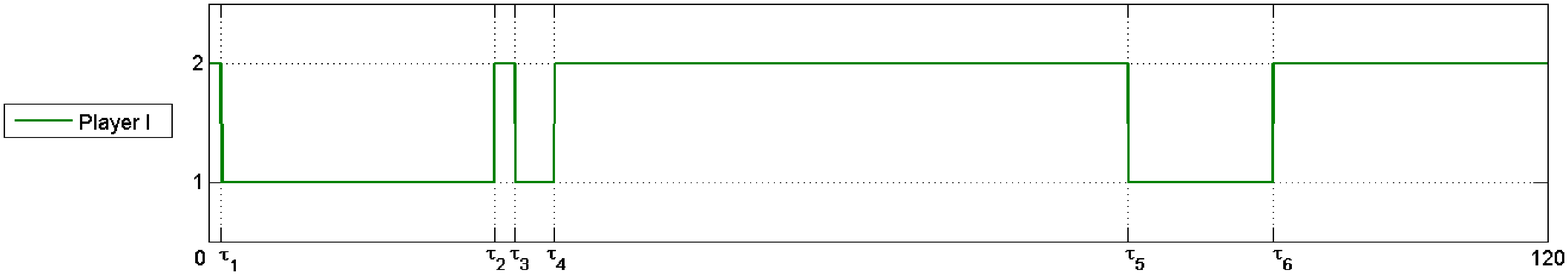}}
{\includegraphics[width=1.02\textwidth]{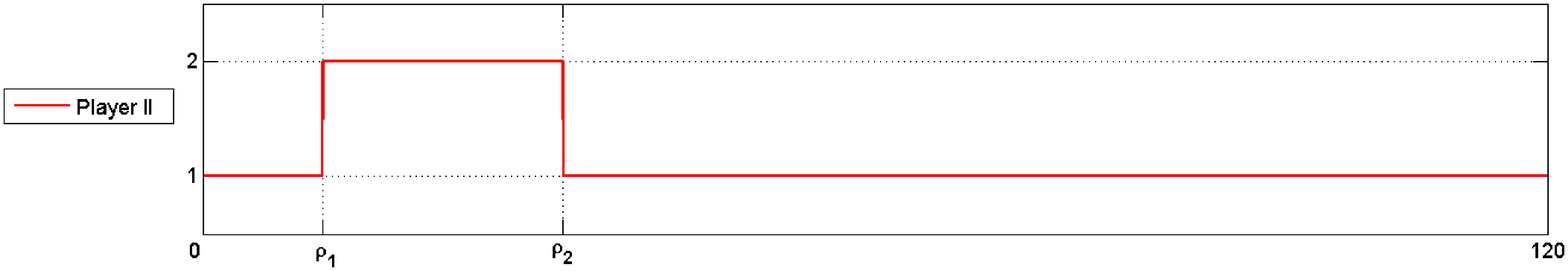}}
\caption{optimal
strategies and state simulation}
\end{figure}
\newpage
\newpage

\end{document}